\documentclass[a4paper,11pt]{article}
\usepackage[hidelinks]{hyperref}
\hypersetup{
	colorlinks=true,   
	linktoc=all,        
	linkcolor=red,     
	citecolor=blue,}
\usepackage{amsmath, amssymb, mathrsfs, calligra}
\usepackage{graphicx}
\usepackage{amssymb}
\usepackage{latexsym, bm}
\usepackage{multicol}
\usepackage{indentfirst}
\usepackage{amssymb,amsfonts}
\usepackage{amsmath}

\usepackage{setspace}
\newtheorem{theorem}{Theorem}[section]
\newtheorem{lemma}{Lemma}[section]
\newtheorem{claim}{Claim}[section]

\newtheorem{corollary}{Corollary}[section]
\newtheorem{proposition}{Proposition}[section]
\newtheorem{definition}{Definition}[section]

\newcommand{\ignore}[1]{}

\begin{document}

\begin{spacing}{0.99}

\title{Ramsey numbers of long even cycles versus books}
\date{}

\author{
{Qizhong Lin\footnote{Center for Discrete Mathematics, Fuzhou University,
Fuzhou, 350108 P.~R.~China. Email: {\tt linqizhong@fzu.edu.cn}. Supported in part  by National Key R\&D Program of China (Grant No. 2023YFA1010202) and NSFC (No.\ 12571361).}}
~~~~and~~~~
Shixi Song\footnote{Center for Discrete Mathematics, Fuzhou University, Fuzhou, 350108 P.~R.~China. Email: {\tt 3503844078@qq.com}.} }

\maketitle

\begin{abstract}
For any positive integers $k$ and $n$, let $B_n^{(k)}$ be the book graph consisting of $n$ copies of the complete graph $K_{k+1}$ sharing a common $K_k$. Let $C_m$ be a cycle of length $m$. Prior work by Allen, \L uczak, Polcyn, and Zhang (2023) established the Ramsey number $R(C_{m},B_n^{(1)})$ for all sufficiently large even integer $m = \Omega(n^{9/10})$. Recently, Hu, Lin, {\L}uczak, Ning, and Peng (2025) obtained the exact value of $R(C_{m},B_n^{(2)})$ under the same asymptotic conditions.
A natural problem is to determine the exact value of $R(C_{m},B_n^{(k)})$ for each fixed $k\ge3$  under similar conditions. This paper provides a complete solution to this problem. 
The lower bound is proved by an explicit construction, while the tight upper bound is established by analyzing the corresponding Ramsey graph using semi-random ideas. 



\end{abstract}

\section{Introduction}\label{sec-1}
Let $\Gamma$ be a graph on $N$ vertices. For any two graphs $G$ and $H$, let $\Gamma\to(G,H)$
denote that either $\Gamma$ contains a copy of $G$ or whose complement contains a copy of $H$. The Ramsey number $R(G,H)$ is the minimum $N$ such that there exists a graph $\Gamma$ on $N$ vertices satisfying $\Gamma\to(G,H)$.  We say $\Gamma$ is a $(G,H)$-Ramsey graph if $\Gamma\nrightarrow(G,H)$.  The existence of the Ramsey number $R(G,H)$ follows from \cite{ram}.

Let $B_n^{(k)}$ be the book graph consisting of $n$ copies of the complete graph $K_{k+1}$ all sharing a common $K_k$.  Books were deeply involved in several recent breakthroughs, such as the exponential improvements \cite{cgms,gnnw} on the diagonal Ramsey numbers, which built upon the book algorithm. 


For a graph $H$ and an integer $N$, the Tur\'{a}n number $ex(N,H)$ is the maximum number of edges in a graph on $N$ vertices that is $H$-free (i.e., it contains no copy of $H$). 
Let $C_m$ be an cycle of length $m$. 
For the even cycle, we only know the correct order of $ex(N,C_{m})$ for $m=4,6,10$, we refer the reader to the survey \cite{fu-s} and references therein. It is widely open of $ex(N,C_{m})$ for each even $m$ of the remaining cases.
In general, for even integer $m$, 
\begin{align*}
N^{1+(4+o(1))/3m}\le ex(N,C_{m})\le O(\sqrt{m}\ln mN^{1+2/m}+m^2N),
\end{align*}
where the upper bound is due to Bukh and Jiang \cite{bj} who improved the classical result due to Bondy and Simonovits \cite{bs}, while
the lower bound by Lubotzky, Phillips, and Sarnak \cite{lps} follows from the construction of regular graphs of large girth.
As observed by Allen, \L uczak, Polcyn, and Zhang \cite{ALPZ-Arxiv} together with a simple induction, we know that  for each fixed $k\ge1$ and even $m\ge4$, 
\begin{align}\label{bou-g}
n+\Omega(n^{(4+o(1))/3m})\le R(C_{m},B_n^{(k)})\le n+O(\sqrt{m}\ln mn^{2/m}+m^2).
\end{align}
There exists a gap between the upper and lower bounds of $R(C_{m},B_n^{(k)})$.

In this paper, we will consider $R(C_m,B_n^{(k)})$, particularly when $m$ is even and sufficiently large. 
For convenience, we adopt the following notation for the remainder of this paper: for each fixed $k \ge 1$, let 
\[
g_k=R(C_m,B_n^{(k)}).
\]

It is known from Dirac \cite{D52} that a graph on $N$ vertices with minimum degree at least $N/2$ is hamiltonian; and if its minimum degree is larger than $N/2$, then it is pancyclic by Bondy \cite{B71}. Therefore, for even $m\ge 2n\ge4$, $g_1=m.$ 
Zhang, Broersma, and Chen \cite{ZBC} determined the value of $g_1$ when $m\ge4$ is even and $3n/4+1\le m\le 2n$. 

Recently,  Allen et al. \cite{ALPZ-Arxiv} established the following result.
\begin{theorem}[Allen, \L uczak, Polcyn, and Zhang \cite{ALPZ-Arxiv}]\label{Thm:alpz}
For every $t \geq 2$, an even integer $m\ge ct^9$ where $c>0$ is a constant, and $n$ satisfying  $(t-1)(m-1) \leq n-1 <t (m-1)$, 
 \[
g_1=\max\{t(m-1)+1,n+\lfloor(n-1)/t\rfloor+1\}.
\]
\end{theorem}

By Theorem \ref{Thm:alpz}, the value of $g_1$ is known for all sufficiently large even integers $$m\ge \Omega(n^{9/10}).$$
On the other hand, as observed in \cite{ALPZ-Arxiv}, Theorem \ref{Thm:alpz} fails if $m\le O(\log n/\log\log n)$ from (\ref{bou-g}).

\ignore{
As stated in \cite{ALPZ-Arxiv}, we may take $m= \Omega(t^9)$ in the above theorem, which holds when $n= O(m^{10/9})$. Therefore, we have the following corollary.
\begin{corollary}[Allen, \L uczak, Polcyn, and Zhang \cite{ALPZ-Arxiv}]
  There exist positive constants $c_1$ and $c_2$ such that for every integer $t \geq 2$, $m\ge c_1 t^9$ and $n\le c_2m^{10/9}$, 
  \[
g_1=\max\{t(m-1)+1,n+\lfloor(n-1)/t\rfloor+1\}.
\]
\end{corollary}
}

For the case $k=2$, we denote this number simply by $R(C_m,B_n)$.
As noted by Faudree, Rousseau, and Sheehan \cite{FRS91}: we know practically nothing  about $R(C_m,B_n)$ when $m$ is even and greater than four. Only recently have we gained some ground on this problem. We refer the reader to \cite{HLLNP,LP21,S10}  and references therein.
In particular, extending Theorem \ref{Thm:alpz}, we obtained the value of $g_2$ for all sufficiently large even  $m\ge\Omega(n^{9/10})$.
\begin{theorem}[Hu, Lin, \L uczak, Ning, and Peng \cite{HLLNP}]\label{hllnp}
For every integer $t \geq 2$, an even integer $m$ sufficiently large, and $n$ satisfying $(t-1)(m-1) \leq n-1 <t (m-1)$,
 \[
g_2=  \begin{cases}
\max\{(t+1)(m-1)+1, n+2\lfloor (n-1)/t \rfloor+2\} &  \textrm{ if } \lfloor (n-1)/t \rfloor+1<m-1,\\
 n+2\lfloor (n-1)/t \rfloor+1 &  \textrm{ if } \lfloor (n-1)/t \rfloor+1=m-1.
\end{cases}
\]
\end{theorem}


In general, the behavior of $g_k$ for $k \ge 3$ would be more difficult to characterize.
We completely determine its value for each $k \ge 3$ under similar conditions,  specifically establishing the value for all sufficiently large even $m$ with $m\ge \Omega(n^{9/10})$.



\begin{theorem}\label{main4}
For all integers $t \geq 2$, $k\ge3$, an even integer $m\ge c(t+k)^9$ where $c>0$ is a constant, and $n$ satisfying $(t-1)(m-1) \leq n-1 <t (m-1)$, the following hold:

  \smallskip
(i) If $\lfloor \frac{n-1}t \rfloor<m-k$, then
$
g_k=\max\{(t+k-1)(m-1)+1, n+k\lfloor (n-1)/t \rfloor+k\}.
$

  \smallskip
(ii) If $\lfloor\frac{n-1}t\rfloor = m-\sigma$ where $ \lceil \frac{k}{2}\rceil+1\le \sigma\le k $, then 
$g_k=n + k \lfloor (n-1)/t \rfloor +\sigma-1.$

  \smallskip
(iii) If $\lfloor\frac{n-1}t\rfloor = m-\sigma$ where $ 2\le \sigma\le \lceil \frac{k}{2}\rceil $, and if
$$g_{k-1}=n+(k-1)\lfloor {(n-1)}/t \rfloor+\ell\;\;\text{where $\sigma-1 \le \ell \le \left\lceil \tfrac{k-1}{2} \right\rceil$,}$$ 
 then, given that $k-1 = \mu_\ell \ell + \alpha_\ell$ with $0\le\alpha_\ell <\ell$, we have 
\[
g_{k}= \begin{cases}
n+k\lfloor (n-1)/t \rfloor+\ell+1 &  \textrm{if $r_k\le r$,}\\
n+k\lfloor (n-1)/t \rfloor+\ell &  \textrm{if  $r_{k}>r$,}
\end{cases}
\]
where $r_{k}=t\lfloor\frac{n-1}t\rfloor+t+k-\ell-n$, and $r=\tfrac{\mu_\ell-1}{\mu_\ell} (t+k-\alpha_\ell)+\alpha_\ell \mu_\ell$.
\end{theorem}

{\em Remark.} 
Upon initial inspection, the seemingly peculiar condition $r_k\le r$ (or $r_{k}>r$) is, in fact, inherent, as can be seen from the remark at the end of the proof of Theorem \ref{main4}.

\medskip
From Theorem \ref{Thm:alpz}, Theorem \ref{hllnp}, and Theorem \ref{main4}(i) for $t=2$, we have the following result.
\begin{corollary}
  For each integer $k\ge1$, there exists a positive integer $n_0=n_0(k)$ such that for all even $n\ge n_0$, $R(C_n,B_n^{(k)})=(k+1)(n-1)+1.$ 
\end{corollary}

\noindent
{\bf Notation.} Let $G=(V,E)$ be a graph. For any subset $U\subseteq V$, we use $G[U]$ to denote the subgraph induced by $U$, and use $G-U$ to denote the subgraph obtained from $G$ by deleting vertices of $U$. For any vertex $v\in V$, let $N_G(v)$ and $d_G(v)$ be the neighborhood and degree of $v$, respectively. Write $N_G(v,U)=N_G(v)\cap U$ and $d_G(v,U)=|N_G(v,U)|$. Let $[n]=\{1,2,\dots,n\}$, and $[m,n]=\{m,m+1,\dots,n\}$.

\medskip\noindent
This paper is organized as follows. In Section \ref{lb}, we will prove lower bounds for Theorem \ref{main4}. In Section \ref{Sec2:Prelimi}, we will collect several previous useful lemmas; and in Section \ref{prop-ex}, we have some properties for $(C_m,B_n^{(k)})$-Ramsey graphs, including a crucial lemma for the proof of the upper bound of Theorem \ref{main4}(iii). Upper bounds for Theorem \ref{main4} will be proven in Section \ref{ub}. In the concluding remarks, a slightly improved bound for $m$ in Theorem \ref{main4} will be presented, together with a discussion of some interesting problems.


\section{Constructions of the lower bounds}\label{lb}

For convenience, let 
\[
p=\lfloor(n-1)/t\rfloor, \;\; \text{and} \;\; n-1=tp+q \;\; \text{where} \;\; 0\le q\le t-1.
\]

(i) Suppose $p<m-k$. We will show that $g_k>\max\{(t+k-1)(m-1), n+kp+k-1\}.$
The construction is similar to that in \cite{ALPZ-Arxiv,HLLNP}.
Let $F$ be the disjoint union of $t+k-1$ cliques of size $m-1$. One can check that $F$ is $C_{m}$-free, and its complement is $B_n^{(k)}$-free since $(t-1)(m-1)<n$ from the assumption. Thus we have $g_k>|V(F)|=(t+k-1)(m-1).$

\begin{figure}[t]
\centerline{
\includegraphics[scale=0.72]{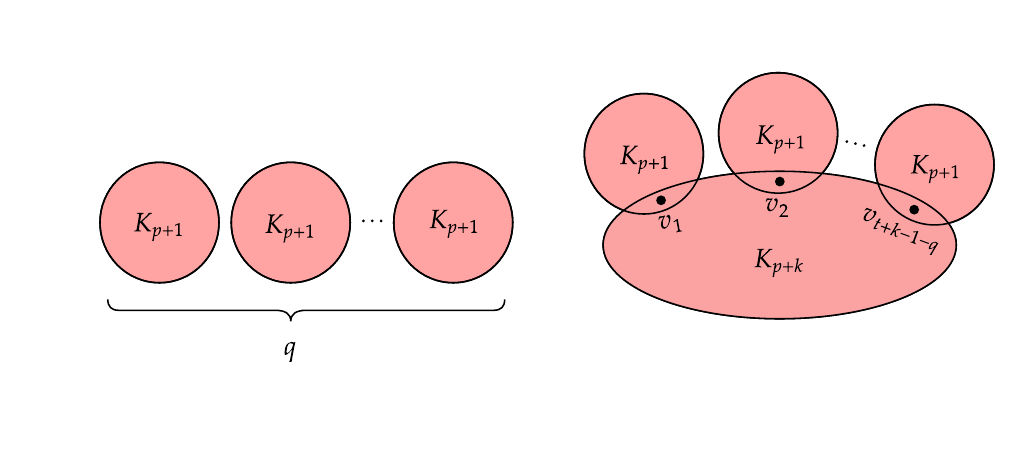}}

\vspace{-1.5cm}
\caption{The graph $F_2$}
\label{fg1}
\end{figure}

It remains to show that $g_k>n+kp+k-1.$
Let $v_1,v_2,\ldots,v_{t+k-1-q}$ be distinct vertices from a $K_{p+k}$.
Let $F_1$ be a graph consisting of the $K_{p+k}$ and extra $t+k-1-q$ vertex-disjoint $K_{p+1}$'s such that the $i$th $K_{p+1}$ shares $v_i$ with the $K_{p+k}$ for $1 \leq i \leq t+k-1-q$. Let $F_{2}$ be a graph obtained from $F_1$ by adding $q$ vertex-disjoint copies of $K_{p+1}$, see Figure \ref{fg1}.
Note that $F_2$ is $C_m$-free as $p+k < m$ from the assumption, and
\begin{align*}
|V(F_2)|=(t+k-1)(p+1)+(p+k)-(t+k-1-q)=n+k(p+1)-1,
\end{align*}
by noting $n-1=tp+q$.

Observe that if $u_1,u_2,\ldots,u_{k}$ form an independent set in $F_2$, then 
$
|\cup_{1\le i\le k}N_{F_2}(u_i)| \geq k p.
$
Thus $u_1,u_2,\ldots,u_{k}$ have at most $|V(F_2)\setminus \cup_{1\le i\le k}N_{F_2}(u_i)|-k=n-1$ common non-neighbors in $F_2$, and so the complement of $F_2$  does not contain $B_n^{(k)}$ as a subgraph. The lower bound follows.

\medskip
(ii) We assume that $p = m - \sigma$, where (for the convenience of proving (iii)) we extend the range of $\sigma$ to $2 \le \sigma \le k$,
and we will show that $g_k>n+kp+\sigma-2.$

Similar to the above, let $v_1,v_2,\ldots,v_{t+k-1-q}$ be distinct vertices from a $K_{p+\sigma-1}$.
Let $F_1$ be a graph consisting of the $K_{p+\sigma-1}$ and extra $t+k-1-q$ vertex-disjoint $K_{p+1}$'s such that the $i$th $K_{p+1}$ shares $v_i$ with the $K_{p+\sigma-1}$ for $1 \leq i \leq t+k-1-q$. Let $F_2$ be a graph obtained from $F_1$ by adding $q$ vertex-disjoint copies of $K_{p+1}$. 
Note that $F_2$ is $C_m$-free as $p+\sigma-1 = m-1$, and
\begin{align*}
|V(F_2)|&=(t+k-1)(p+1)+(p+\sigma-1)-(t+k-1-q)\\
&=(t+k)p+\sigma-1+q\\
&=n+kp+\sigma-2.
\end{align*}
Observe that if $u_1,u_2,\ldots,u_{k}$ form an independent set in $F_2$, then 
\[
|\cup_{1\le i\le k}N_{F_2}(u_i)| \geq(k-1)p+(p+\sigma-2)-(k-1)= k (p-1)+\sigma-1,
\]
where the minimum will be achieved if some vertex $u_i$ was chosen from the $K_{p+\sigma-1}$ other than $v_1,v_2,\ldots,v_{t+k-1-q}$.
Therefore, $u_1,u_2,\ldots,u_{k}$ have at most $n-1$ common non-neighbors in $F_2$, and so the complement of $F_2$ contains no $B_n^{(k)}$. Thus we obtain $g_k > n+kp+\sigma-2$.

\medskip
(iii) For this case, we assume that 
\[
p=m-\sigma,\;\;\text{ where $2\le \sigma\le \lceil\tfrac{k}{2}\rceil$.}
\]
Furthermore, we assume that for $k\ge3$,
$ g_{k-1}\ge n+(k-1)p+\ell, \;\text{where}\; \sigma-1\le \ell\le \lceil\tfrac{k-1}{2}\rceil.$ Given that 
$k-1=\mu_\ell \ell+\alpha_\ell$ where $0\le \alpha_\ell<\ell$, we will show that
\begin{align}\label{gk}
 g_{k}\ge \begin{cases}
n+kp+\ell+1 &  \textrm{if $r_{k}\le r$,}\\
n+kp+\ell &  \textrm{if  $r_{k}>r$,}
\end{cases}
\end{align}
where $r_{k}=tp+t+k-\ell-n$, and $r=\tfrac{\mu_\ell-1}{\mu_\ell} (t+k-\alpha_\ell)+\alpha_\ell \mu_\ell$.

 
First we will show (\ref{gk}) holds for $k=3$. For this case, we know $\sigma=2$, and so $p+1=m-1$. By Theorem \ref{hllnp},  $g_2=n+2p+1$. Thus $\ell=1$, $\mu_1=2$ and $\alpha_1=0$. Also, $r_3=tp+t+2-n$.

Assume $r_3\leq\tfrac{t+3}{2}$.
We will prove $g_3> n+3p+1$.
Let $\Lambda$ be the graph consisting of $2r_3$ copies of $K_{p+1}$'s paired up in twos where each such pair of $K_{p+1}$'s shares a common vertex, and let $\Gamma_3$ be obtained from $\Lambda$ by adding extra $t+3-2r_3$ vertex-disjoint $K_{p+1}$'s. 
Note that $\Gamma_3$ is $C_m$-free as $p+1= m-1$ from the assumption. Moreover,  
$
|V(\Gamma_3)|=(t+3)(p+1)-r_3=n+3p+1.
$
Observe that if $u_1$, $u_2$ and $u_3$ form an independent set in $\Gamma_3$, then 
$|\cup_{1\le i\le 3}N_{\Gamma_3}(u_i)| \geq 3p-1.$
Therefore, $u_1$, $u_2$ and $u_3$ have at most $n-1$ common neighbors in
the complement of $\Gamma_3$, and so the complement of $\Gamma_3$ is $B_n^{(3)}$-free.
The lower bound follows for this case.

Now assume $r_3>\tfrac{t+3}{2}$, and we will show that $g_3> n+3p$.
Let $\Gamma_3'$ be the graph obtained from $q$ vertex-disjoint $K_{p+1}$'s by adding extra $t+3-q$ copies of $K_{p+1}$'s with a common vertex.
Clearly $\Gamma_3'$ is $C_m$-free as $p+1=m-1$.
Moreover, recall $q=n-1-tp$, so we have that
\begin{align*}
|V(\Gamma_3')|=(t+3)(p+1)-(t+2-q)=tp+q+3(p+1)-2=n+3p.
\end{align*}
Note that if $u_1$, $u_2$ and $u_3$ form an independent set in $\Gamma_3'$, then
$
|\cup_{1\le i\le3} N_{\Gamma_3'}(u_i)| \geq 3p-2.
$
Therefore, $u_1$, $u_2$ and $u_3$ have at most $n-1$ common neighbors in
the complement of $\Gamma_3'$, and so the complement of $\Gamma_3'$ is $B_n^{(3)}$-free.
So the assertion holds for $k=3$.

\begin{figure}[t]
\begin{center}
    {\includegraphics[scale=0.8]{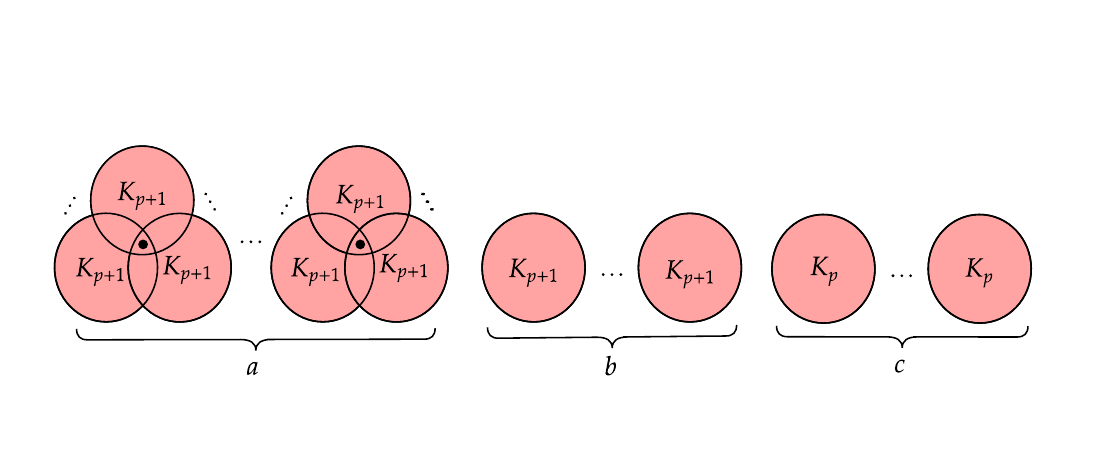}}
\end{center}

\vspace{-1.5cm}
\caption{The graph $\Gamma_{a,b,c}$}
\label{fg6}
\end{figure}

\medskip
In the following, we assume $k\ge4$ and suppose that (\ref{gk}) holds for smaller $k$. We will show that (\ref{gk}) also holds for $k$.
We first define a family of graphs $\Gamma_{a,b,c}$, see Figure \ref{fg6}.
\begin{definition} 
Let $a, b, c$ be non-negative integers. The graph $\Gamma_{a,b,c}$ is constructed as follows:

\smallskip
    (1) It contains $a$ clusters, where each cluster is formed by some copies of $K_{p+1}$ with one common vertex.

    (2) It has $b$ pairwise vertex-disjoint copies of $K_{p+1}$ not intersecting any cluster.

    (3) It has $c$ pairwise vertex-disjoint copies of $K_{p}$ not intersecting any cluster or $K_{p+1}$.

\end{definition}

\ignore{
For the base case $k=2$, we establish the validity of the proof along with a type $\Gamma_{a,b,c}$ lower-bound construction. Under the inductive hypothesis that the proof remains valid for a general $k$ with a corresponding type $\Gamma_{a,b,c}$ construction, we proceed to analyze the case for $k+1$.}


\noindent
{\bf Case 1.} \ $r_{k}\le r$.

\vspace{0.2cm}

\begin{figure}[t]
\centerline{
\includegraphics[scale=0.85]{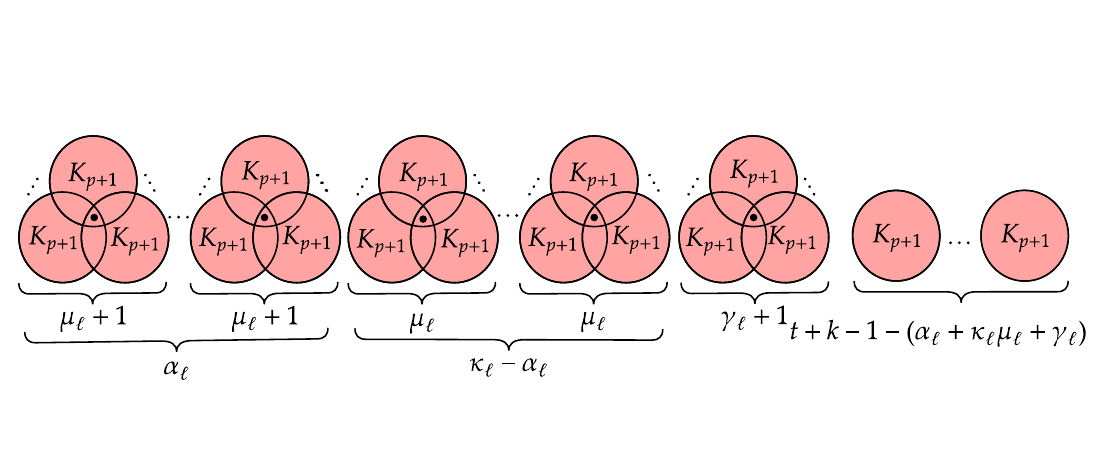}}

\vspace{-1cm}
\caption{The graph $\Gamma_{k}$}
\label{fg7}
\end{figure}

Observe that $r_{k}=tp+t+k-\ell-n\ge k-\ell$. Moreover, we have $k-1=\mu_\ell \ell+\alpha_\ell$ with $0\le \alpha_\ell<\ell$, and $\ell\le\lceil\frac{k-1}{2}\rceil$. So we can express $r_{k}-\alpha_\ell$ as $(\mu_\ell-1)\kappa_\ell+\gamma_\ell$, where $\kappa_\ell\ge\ell> \alpha_\ell$ and $0\le \gamma_\ell\leq \mu_\ell-2$.
Let $G_1$ be a graph consisting of $\mu_\ell+1$ copies of $K_{p+1}$ with a common vertex,  $G_2$ a graph consisting of $\mu_\ell$ copies of $K_{p+1}$ with a common vertex, and $G_3$ a graph consisting of $\gamma_\ell+1$ copies of $K_{p+1}$ with a common vertex. Let $\Gamma_k$ be the union graph consisting of $\alpha_\ell$ disjoint copies of $ G_1$, $(\kappa_\ell-\alpha_\ell)$ disjoint copies of $G_2$, one copy of $ G_3$, and extra $t+k-1-(\alpha_\ell+\kappa_\ell \mu_\ell+\gamma_\ell)$ copies of $K_{p+1}$, and all of these $G_1$'s, $G_2$'s, $G_3$'s, and $K_{p+1}$'s are vertex-disjoint from each other, see Figure \ref{fg7}.
Clearly, $\Gamma_{k}$ is $C_m$-free since $p+1 \le m-1$. We will verify that the complement of $\Gamma_{k}$ is $B_n^{(k)}$-free.
Note that the total number of vertices in these $G_i$'s for $1\le i\le3$ is counted
\begin{align*}
&\alpha_\ell((\mu_\ell+1)(p+1)-\mu_\ell)+(\kappa_\ell-\alpha_\ell)(\mu_\ell (p+1)-(\mu_\ell-1))+(\gamma_\ell+1)(p+1)-\gamma_\ell\\
         &=(p+1)(\alpha_\ell(\mu_\ell+1)+(\kappa_\ell-\alpha_\ell)\mu_\ell+\gamma_\ell+1)-\alpha_\ell\mu_\ell-(\kappa_\ell-\alpha_\ell)(\mu_\ell-1)-\gamma_\ell\\
         &=(p+1)(\alpha_\ell+\kappa_\ell\mu_\ell+\gamma_\ell+1)-(\alpha_\ell+\kappa_\ell(\mu_\ell-1)+\gamma_\ell)\\
       &=(p+1)(\alpha_\ell+\kappa_\ell\mu_\ell+\gamma_\ell+1)-r_{k}.
\end{align*}
Thus $|V(\Gamma_{k})|=(t+k)(p+1)-r_{k}=k(p+1)+(n-k+\ell)=n+kp+\ell$ by noting $r_{k}=tp+t+k-\ell-n$.
We aim to compute, for an independent set ${u_1, \ldots, u_k}$ in $\Gamma_k$, the minimum possible size of $| \cup_{i=1}^k N_{\Gamma_k}(u_i) |$. This minimum is achieved in the worst case, when none of the $u_i$ is a center vertex of any $G_j$ ($1 \le j \le 3$). Recall $k-1=\mu_\ell\ell+\alpha_\ell$ with $0\le \alpha_\ell<\ell$.  Therefore,
\begin{align*}
|\cup_{1\le i\le k}N_{\Gamma_{k}}(v_{i})| &\geq\alpha_\ell((\mu_\ell+1)p+1)+(\ell-\alpha_\ell)(\mu_\ell p+1)+p+1-k
\\&=(\mu_\ell\ell+\alpha_\ell)p+\ell+p+1-k
\\&=k(p-1)+\ell+1.
\end{align*}
 It follows that $u_1,u_2,\ldots,u_{k}$  have at most $n-1$ common neighbors in
the complement of $\Gamma_{k}$, and so the complement of $\Gamma_{k}$ is $B_n^{(k)}$-free. 


\medskip
\noindent
{\bf Case 2.} \ $r_{k}>r$.

\vspace{0.2cm}
We aim to show $g_k>n+kp+\ell-1$. If $\ell=\sigma-1$, then we can apply $F_2$ as in (ii) (which is not of type $\Gamma_{a,b,c}$) to obtain a lower bound for $g_{k}$, i.e., $g_{k}> n+kp+\sigma-2=n+kp+\ell-1$.
So we may assume that $\sigma\le\ell\le\lceil\tfrac{k-1}{2}\rceil$. 
\begin{claim}\label{construction} 
If $\Gamma_{k-1}$ is of type $\Gamma_{a,b,c}$ on $n+(k-1)p+\ell-1$ vertices and $\Gamma_{k-1}\nrightarrow(C_m,B_n^{(k-1)})$,
and $\Gamma_{k}$ is obtained from $\Gamma_{k-1}$ by adding a $K_{p}$, then $\Gamma_{k}$ is of type $\Gamma_{a,b,c}$ with $n+kp+\ell-1$ vertices satisfying $\Gamma_{k}\nrightarrow(C_m,B_n^{(k)})$.
\end{claim}
\noindent
{\bf Proof.}
From the construction, $\Gamma_{k}$ is of type $\Gamma_{a,b,c}$ with $n+kp+\ell-1$ vertices. Moreover, $\Gamma_{k}$  is $C_m$-free. It remains to show that the complement of $\Gamma_{k}$ is $B_n^{(k)}$-free.
Since there is no $B_n^{(k-1)}$ in the complement of $\Gamma_{k-1}$, we have that for any $k-1$ vertices $u_1,u_2,\ldots,u_{k-1}$ that form an independent set in $\Gamma_{k-1}$, 
$
|\cup_{1\le i\le k-1} N_{\Gamma_{k-1}}(u_i)| \geq |V(\Gamma_{k-1})|-(k-1)-(n-1)=(k-1)(p-1)+\ell.
$
Since $\Gamma_{k}$ is of type $\Gamma_{a,b,c}$, any vertices $u_1,u_2,\ldots,u_{k}$ that form an independent set in $\Gamma_{k}$ must satisfy that 
$
 |\cup_{1\le i\le k} N_{\Gamma_{k}}(v_i)| \geq (k-1) (p-1)+\ell+(p-1)=k(p-1)+\ell.
$
Therefore, $u_1,u_2,\ldots,u_{k}$ have at most $n-1$ common neighbors, and so the complement of $\Gamma_{k}$ is $B_n^{(k)}$-free as desired. \hfill $\square$

\medskip
Note that we have proved (\ref{gk}) for $k=3$. Specifically, $\sigma=2$ when $k=3$. For $r_3>\frac{t+3}{2}$, we know $g_3>n+3p$, which corresponds to $\ell=\sigma-1=1$. For $r_3\le\frac{t+3}{2}$, there exists a $\Gamma_{3}$ of type $\Gamma_{a,b,c}$ on $n+3p+1$ vertices such that $\Gamma_{3}\nrightarrow(C_m,B_n^{(3)})$. In the following, suppose $k\ge 4$ and that for $2\le \sigma\le \lceil \frac{k-1}{2}\rceil$, there exists a $\Gamma_{k-1}$ of type $\Gamma_{a,b,c}$ on $n+(k-1)p+\ell-1$ vertices (except when $\ell=\sigma-1$) with $\Gamma_{k-1}\nrightarrow(C_m,B_n^{(k-1)})$. We will show this holds for $k$.

If $k$ is even, then $2\le\sigma\le \lceil \frac{k}{2}\rceil=\lceil \frac{k-1}{2}\rceil$. 
Moreover, if $r_k\le r$, then by Case 1, there exists a $\Gamma_{k}$ of type $\Gamma_{a,b,c}$ on $n+kp+{\ell}$ vertices such that $\Gamma_{k}\nrightarrow(C_m,B_n^{(k)})$.  If $r_k>r$, by the induction hypothesis and Claim \ref{construction}, there exists a $\Gamma_{k}$ of type $\Gamma_{a,b,c}$ on $n+kp+\ell-1$ vertices  (except when $\ell=\sigma-1$) that satisfies $\Gamma_{k}\nrightarrow(C_m,B_n^{(k)})$. 

If $k$ is odd, we first consider the case where $\sigma=\lceil\tfrac{k}{2}\rceil=\tfrac{k+1}{2}=\lceil\tfrac{k-1}{2}\rceil+1$. By the lower bound of Theorem \ref{main4}(ii), we know $g_{k-1}\ge n+(k-1)p+(\sigma-1)$. We may assume $\ell=\sigma-1$, as assuming a larger value for $\ell$ would only strengthen the lower bound. Now if $r_k\le r $, then applying the result from Case 1 we get a $\Gamma_{k}$ of type $\Gamma_{a,b,c}$  on $n+kp+(\sigma-1)$ vertices satisfying $\Gamma_{k}\nrightarrow(C_m,B_n^{(k)})$. If $r_k>r$, then as observed at the beginning of Case 2, we have $g_{k}>n+kp+\sigma-2$. 

Now we may assume $2\le \sigma\le\lceil\tfrac{k-1}{2}\rceil$. 
 If {$r_k\le r$,} then by Case 1 there exists a $\Gamma_{k}$ of type $\Gamma_{a,b,c}$ with $n+kp+\ell$ vertices satisfying that $\Gamma_{k}\nrightarrow(C_m,B_n^{(k)})$.  If {$r_k>r$,} then by induction and Claim \ref{construction} there exists a $\Gamma_{k}$ of type $\Gamma_{a,b,c}$ on $n+kp+\ell-1$ vertices (except when $\ell=\sigma-1$)  satisfying $\Gamma_{k}\nrightarrow(C_m,B_n^{(k)})$. 
 
 We complete the proof of the lower bounds of Theorem \ref{main4}. \hfill$\Box$


\section{Useful lemmas}\label{Sec2:Prelimi}
In the following, we will focus on the upper bounds of Theorem \ref{main4}. The proofs proceed by induction on $k\ge2$, analyzing the $(C_m,B_n^{(k)})$-Ramsey graph where $m$ is even, i.e., $R$ contains no even cycle $C_m$ and whose complement is $B_n^{(k)}$-free. The strategy is twofold. First, if the graph has sufficiently large minimum degree, it must contain consecutive even cycles across a wide range, ensuring the existence of a $C_m$. This approach has been successfully applied in \cite{ALPZ-Arxiv,HLLNP,LP21} and even to a problem in spectral graph theory \cite{LN-Arxiv}. Second, we utilize the extremal properties of such Ramsey graphs (see Section \ref{prop-ex}) to derive a contradiction when the graph's order is large enough.

 A graph $G$ is Hamiltonian  if it contains a cycle that visits every vertex exactly once; such a cycle is called a Hamiltonian cycle. Let $\delta(G)$ be the minimum degree in $G$, and let $c(G)$ and $ec(G)$ denote the length of a longest cycle and a longest even cycle in $G$, respectively. We say that a graph $G$ is pancyclic if it contains cycles of every length from 3 to $|V(G)|$. We begin by reviewing the following classical results.

\begin{lemma}[Dirac \cite{D52}] \label{dirac}
If an $n$-vertex graph $G$ is 2-connected, then $c(G)\geq \min\{2\delta(G),n\}$.
\end{lemma}

Bondy obtained a stronger result if graph $G$ has minimum degree $\delta(G) \geq |V(G)|/2$.
\begin{lemma} [Bondy \cite{B71}]\label{pancyclic}
If an $n$-vertex graph $G$ has minimum degree $\delta(G) \geq n/2$, then either $G$ is pancyclic or $G$ is the complete bipartite graph $K_{n/2,n/2}$.
\end{lemma}

The following result guarantees long even cycles in graphs.

\begin{lemma}[Voss and Zuluaga \cite{VZ77}]\label{VZ}
If $G$ is a 2-connected graph on $n$ vertices, then we have $ec(G) \geq \min\{2\delta(G),n-1\}$.
\end{lemma}

 We also need the following result on the existence of consecutive even cycles.
\begin{lemma}[Gould, Haxell, and Scott \cite{GHS02}]\label{GHS}
Let $\varepsilon> 0$, $K=75\cdot 10^4/\varepsilon^5$, and let $G$ be a graph with $n\geq \frac{45K}{\varepsilon^4}$
vertices and minimum degree at least $\varepsilon n$. Then $G$ contains a cycle of length $t$ for
every even integer $t\in [4,ec(G)-K]$.
\end{lemma}

The following result shows that if the minimum degree of a graph is large enough, then for any pair of vertices there are paths of consecutive lengths.
\begin{lemma}[Williamson \cite{will}] \label{will}
Let $G$ be a graph with $n$ vertices. If $\delta(G) \geq \tfrac{n}{2}+1$,
then for any two vertices $x,y$ and any $2 \leq \ell \leq n-1$, there is an $(x,y)$-path of length $\ell$  in $G$.
\end{lemma}

The following lemma states that if the minimum degree of a graph $G$ is large enough, then we can find a small set $W$ such that $G-W$ is a union of some vertex-disjoint 2-connected subgraphs.
\begin{lemma}[Allen,
\L uczak, Polcyn, and Zhang \cite{ALPZ-Arxiv}] \label{Luzack}
Let $n \geq \ell \geq 2$. For each graph $G$ with $n$ vertices and
minimum degree $\delta(G) \geq n/\ell+\ell$,  there exists an $s<\ell$
and a subset $W \subset V(G)$ with $|W| \leq s-1$, such that $G-W$ is a union
of  $s$ vertex-disjoint 2-connected graphs.
\end{lemma}

\section{Properties for the Ramsey graphs}\label{prop-ex}

In this section, we shall have some properties for the $(C_m,B_n^{(k)})$-Ramsey graph $R$. 
For any $t\ge 2$ and $k\ge3$, let  $\varepsilon=\tfrac{1}{2(t+k)}$, and let
\[
K=2.4\cdot 10^7 (t+k)^5, \;\;\text{and} \;\; \text{even}\;\;m\geq \frac{45K}{\varepsilon^4}=1.728\cdot 10^{10}(t+k)^9.
\]

{\bf Assumption:} Let $s\ge k$ be an integer.
Assume that there exists a set $W\subset V(R)$ with $|W|\le s-1$ such that $R-W$ is a union
of  $s$ vertex-disjoint 2-connected graphs induced by $V_1,\ldots,V_s$ respectively. Then $W,V_1,\ldots,V_s$ form a partition of $V(R)$.
For each $1 \leq i \leq s$, let 
$$W_i=\{w \in W: d_{R}(w,V_i) \geq 2s\}, \;\; \text{and} \;\; U_i=V_i \cup W_i.$$ 
Let $F$ be a bipartite graph with vertex set $V(F)=W \cup \{U_1,\ldots, U_s\}$ and $\{w,U_i\} \in E(F)$ if and only if $w \in U_i$. 

The following lemma can be found in \cite{ALPZ-Arxiv,HLLNP}; we include a proof for completeness. 

\begin{lemma}\label{forest}
If $|V_i \cup V_j| \geq m$ for any $1\leq i \neq j \leq k$, and any pair of vertices of $V_i$ for $1 \leq i \leq k$ is connected by a path of length from two to $m/2$ using vertices from $V_i$, then $F$ is a forest.
\end{lemma}

\noindent
{\bf Proof.} Suppose to the contrary that  $w_{i_1},U_{j_1},w_{i_2},U_{j_2},\ldots,w_{i_r},U_{j_r},w_{i_1}$ is a cycle of length $2r$ in $F$. Since for each $1 \leq \ell \leq r$, $w_{i_\ell}$ has at least two neighbors in  $V_{j_\ell}$,  there are distinct vertices $x_\ell,y_\ell \in V_{j_\ell}$ such that $x_\ell$ is a neighbor of $w_{i_\ell}$ and $y_\ell$ is a neighbor of $w_{i_{\ell+1}}$. Here, $w_{i_{r+1}}=w_{i_1}$. Let $P_\ell$ be a path connecting $x_\ell$ and $y_\ell$ using vertices from $V_{i_\ell}$. Thus  $w_{i_1},x_1,P_1,y_1,w_{i_2},x_2,P_2,y_2,\ldots,w_{i_r},x_r,P_r,y_r,w_{i_1}$ is a cycle in $R$. By the assumption, the length of $P_\ell$ can vary from two to $m/2$. So we can get a $C_{m}$ by choosing $P_\ell$ properly, a contradiction. Therefore, $F$ is a forest. \hfill $\square$

\medskip
An analogous argument to \cite{ALPZ-Arxiv,HLLNP} yields the following result under a slightly relaxed condition.
\begin{lemma}\label{basic}
 Suppose the minimum degree $\delta(R) > m/2+K/2+s$, where $K$ is the parameter derived from Lemma \ref{GHS}. We have the following properties.
  
  \smallskip
  (i) $V=\cup_{i=1}^s U_i$.
  
  \smallskip
  (ii) $|U_i| \leq m-1$ for each $1 \leq i \leq s$.
  
  \smallskip
  (iii) $\sum_{i=1}^{s}|U_i| \leq |V|+s-1$.
  
\end{lemma}

\noindent
{\bf Proof.} 
(i) It suffices to show that each $w \in W$ is contained in some $U_i$. Otherwise, there exists a vertex $w\in W$ that is not in $U_i$ for each $1 \leq i \leq s$. Then $d_R(w) < 2s^2+|W|\le2s^2+s<m/2+K/2+s$ provided $m\ge 4s^2$, a contradiction. 

\smallskip
(ii) We first show $|V_i| \leq m-1$. From the assumption, $\delta(R[V_i]) \geq \delta(R)-|W| > m/2+K/2$. 
If $|V_i| \le2\delta(R[V_i])$, then Lemma \ref{pancyclic} implies that there is a $C_m$.
Else if $|V_i| > 2\delta(R[V_i])$, then by Lemma \ref{VZ}, there is an even cycle of length at least $2\delta(R[V_i])\ge m+K$ in $R[V_i]$.
Thus $C_m$ is guaranteed by applying Lemma \ref{GHS} with $\varepsilon=\tfrac{1}{2(t+k)}$, a contradiction. 

In the following, we will show $|U_i| \leq m-1$ for $1 \leq i \leq s$.
Otherwise, suppose $|U_1| \geq m$ without loss of generality. Take a subset $U_1' \subset U_1$ with $V_1 \subset U_1'$ and $|U_1'|=m$. Let $H=R[U_1']$ and we next show $H$ is  Hamiltonian. The fact $U_1' \setminus V_1 \subset W$ implies $|U_1'\setminus V_1| <s$, and so $|V_1| > m-s$.
 By joining a pair of non-adjacent vertices with degree sum (in $H$) at least $m$ iteratively, we get the $m$-closure of $H$, denoted by $H'$.  Recall $\delta(R[V_1]) \geq \delta(R)-s >m/2$. It follows that $V_1$ forms a clique in $H'$.
Assume $U_1'\setminus V_1=\{w_1,\ldots,w_\ell\}$ with $\ell < s$. From the assumption $d_R(w,V_1) \geq 2s$ for each $w \in U_1'\setminus V_1$, we are able to find $x_i,y_i \in V_1$ for $1 \leq i \leq \ell$ such that $\{x_i,y_i\} \subset N_R(w_i,V_1)$ and they are all distinct. As $V_1$ forms a clique in $H'$, we can find a $C_m$ in $H'$.
 Therefore, $H$ is Hamiltonian, which implies that $U_1$ contains $C_m$, a contradiction. 
  
\smallskip
(iii) Note that each vertex from $V_i$ has at least $\delta(R)-s>m/2$ neighbors in $V_i$. Thus $|V_i| >m/2$ and $|V_i \cup V_j| \geq  m$. From (ii), $|V_i| \leq |U_i| \leq m-1$, we get $\delta(R[V_i]) \geq \tfrac{|V_i|}{2}+1$. It follows by Lemma \ref{will} that each pair of vertices in $V_i$ is connected by a path of length from $2$ to $m/2$. Thus, by Claim \ref{forest}, $F$ is a forest. 

In the sum $\sum_{i=1}^{s}|U_i|$, each vertex $w \in W$ is counted exactly $d_F(w)$ times and each vertex in $\cup_{i=1}^s V_i=V \setminus W$ is counted exactly once. Thus,
$\sum_{i=1}^{s}|U_i|=|V|-|W|+\sum_{w \in W} d_F(w).$
Note that $F$ is a forest, we have $\sum_{w \in W} d_F(w) \leq |W|+s-1$. Therefore, $\sum_{i=1}^{s}|U_i| \leq |V|+s-1$. \hfill $\square$

\begin{lemma}\label{c5}
If $\delta(R)\ge 2s^2+s$, then for any distinct $i_1,\dots,i_k\in[s]$,  $|V \setminus(\cup_{1\le j\le k} U_{i_j})| \leq n-1$.
\end{lemma}
{\bf Proof.} For each $w \in W$ with $w \not \in W_{i_1}$, it has at most $2s$ neighbors in $V_{i_1}$. Thus vertices from $W\setminus W_{i_1}$ have at most $2s^2$ neighbors in $V_{i_1}$. As $|V_{i_1}| \geq \delta(R)-|W|>2s^2$, there is a vertex $v_{i_1}\in V_{i_1}$ such that $v_{i_1}$ is non-adjacent to vertices in $W \setminus W_{i_1}$. Similarly, we can find $v_{i_2} \in V_{i_2},  \dots, v_{i_k} \in V_{i_k}$ such that $v_{i_j}$ is non-adjacent to vertices in $W \setminus W_{i_j}$ for $2\le j\le k$. Therefore, $v_{i_1}, v_{i_2}, \dots, v_{i_k}$ are non-adjacent to all vertices in $W\setminus (\cup_{1\le j\le k} W_{i_j})$ and so in $V\setminus (\cup_{1\le j\le k} U_{i_j})$.
As we assumed that the complement of $R$ is $B_n^{(k)}$-free, we obtain $|V \setminus(\cup_{1\le j\le k} U_{i_j})| \leq n-1$ as claimed.\hfill$\square$

\begin{definition}[Duplication number]
For any $\mathscr{U}=\{U_1,U_2,\dots,U_h\}$, we call the difference $\sum_{i=1}^{h}| U_{i}|-|\cup_{1\le i\le h} U_{i}|$ the duplication number of $\mathscr{U}$. Also, we call $\frac{1}{h}(\sum_{i=1}^{h}| U_{i}|-|\cup_{1\le i\le h} U_{i}|)$ the duplication rate of $\mathscr{U}$. 

\end{definition}

The reader is advised to note carefully the above definitions.

\begin{lemma}\label{dp-u}
For any $h$ sets $U_1,U_2,\dots,U_h$, if  $|U_i\cap (\cup_{1\le j\neq i\le h}U_j)|\ge 2$ for each $1\le i\le h$, then $\sum_{1\le i\le h}|U_i|-|\cup_{1\le i\le h}U_i|\ge h$.
\end{lemma}
{\bf Proof.} Since $|U_i \cap (\cup_{1\le j \neq i\le h} U_j)| \ge 2$ for each $1 \le i \le h$, each $U_i$ contains at least two vertices that are also in other sets.
Let $x_i$ be the number of vertices in $U_i$ that belong to $\cup_{1\le j\neq i\le h}U_j$.
Let $x = \sum_{i=1}^h x_i$. Then $x \ge 2h$ because $x_i \ge 2$ for each $i$.
Suppose there are $y$ vertices that are in at least two sets. Then $y \le x/2$ because each such vertex is counted in at least two of the $x_i$'s. Then $\sum_{i=1}^h |U_i| - |\cup_{i=1}^h U_i| = x - y \ge x - x/2 = x/2 \ge h$.\hfill$\Box$

\medskip
The following lemma is crucial for the proof of the upper bound of Theorem \ref{main4}(iii), which states that if any $k$ elements of a family of sets have small duplication number, then the size of the union of these sets from the family would be large.

\begin{lemma}\label{c6}
Suppose that $|U_i|=p+1$ for each $i\in[ t+k]$. For $ \ell\ge 1$, suppose $k-1 = \mu \ell + \alpha$ where $0\le\alpha <\ell$.
If $|\cup_{1\le j\le k} U_{i_j}|\geq kp+\ell+1$ for any distinct $i_1,i_2,\dots,i_k\in[ t+k]$, then $|V|\geq \sum_{i=1}^{t+k}| U_i|-\tfrac{\mu-1}{\mu}(t+k-\alpha)-\alpha\mu$.
\end{lemma}
{\bf Proof.} Let
$
\mathscr{U}=\{U_1,U_2,\dots, U_{t+k}\}.
$
As $|U_i|=p+1$ for $i\in[t+k]$ and the union of any $k$ elements of $\mathscr{U}$ is at least $kp+\ell+1$, we have that the duplication number of any $k$ elements of $\mathscr{U}$ is at most $k-\ell-1$. 

\ignore{Suppose there exist $k-1$ elements of $\mathscr{U}$, say $U_1,\dots,U_{k-1}$, such that their duplication number is $k-\ell-1$. Then for any $i\in[k-1]$ and $j\in[k,t+k]$, $U_i\cap U_j=\varnothing$. Otherwise, $U_1,\dots,U_{k-1}$ together with some $U_j$ for $j\in[k,t+k]$ have the duplication number at least $k-\ell$, a contradiction. 

Let $b$ be the maximum duplication number of any $t+1$ elements of $\mathscr{U}$ that may achieve. Then 
\begin{align}\label{max-dup}
\sum_{i=1}^{t+k}| U_{i}|-|V|\le  k-\ell-1+b.
\end{align}
}

\begin{claim}\label{max-dup-1}
Let $b$ be the maximum duplication number of any $t+1$ elements of $\mathscr{U}$ that may achieve. Then $\sum_{i=1}^{t+k}| U_{i}|-|V|\le k-\ell-1+b$.
\end{claim}

\noindent
{\bf Proof.} Suppose first that there exist $k-1$ elements of $\mathscr{U}$, say $U_1,\dots,U_{k-1}$, such that their duplication number is $k-\ell-1$. Then for any $i\in[k-1]$ and $j\in[k,t+k]$, $U_i\cap U_j=\varnothing$; otherwise, $U_1,\dots,U_{k-1}$ together with some $U_j$ for $j\in[k,t+k]$ have the duplication number at least $k-\ell$, a contradiction. So $\sum_{i=1}^{t+k}| U_{i}|-|V|\le  (k-\ell-1)+b$, the assertion follows for this case. 

Thus, we only need to verify that the assertion  holds provided that any $k-1$ elements of $\mathscr{U}$ have duplication number at most $k-\ell-2$. There are two subcases.

\medskip
(i) There exist $k-2$ elements of $\mathscr{U}$, say $U_1,\dots,U_{k-2}$, such that their duplication number is $k-\ell-2$. 
Let  $\mathscr{U}_{k-2}=\{U_{k-1},\dots,U_{t+k}\}$ for convenience.

\begin{proposition}\label{indt}
$\sum_{i=k-1}^{t+k}| U_{i}|-|\cup_{k-1\le i\le t+k}U_i|\leq b+1$.     
\end{proposition}
{\bf Proof.} 
We first claim that for some $U_i$ with $i\in[k-1, t+k]$, $|U_i\cap (\cup_{k-1\le j\neq i\le t+k}U_j)|\le 1$. Otherwise,  for each $i\in[k-1, t+k]$, $|U_i\cap (\cup_{k-1\le j\neq i\le t+k}U_j)|\ge 2$. By Lemma \ref{dp-u}, the duplication number of $\mathscr{U}_{k-2}$ is at least $t+2$.

If $t+2\geq k$, then we can delete from $\mathscr{U}_{k-2}$, one by one, those subsets that have the fewest common vertices with the remaining subsets. We delete a subset from $\mathscr{U}_{k-2}$ that has the fewest common elements with the remaining subsets. Then the remaining $t+1$ subsets of $\mathscr{U}_{k-2}$ have at least $t$ duplication numbers. If the $t+1$ subsets have exactly $t$ duplication numbers, then there exists a subset that shares at most one vertex with the others; deleting it leaves $t$ subsets with at least $t-1$ duplication numbers.
If instead the $t+1$ subsets have at least $t+1$ duplication numbers, then we proceed similarly to the initial case to obtain $t$ subsets with at least $t-1$ duplication numbers.
We repeat this process until we have deleted subsets one by one down to a certain number. Finally, we obtain $k$ sets $U_{i_1}, U_{i_2}, \dots,U_{i_k}$ for $i_j\in[k-1, t+k]$ whose duplication number is at least
 $k-1>k-\ell-1,$ which contradicts our assumption.
 
So we may assume that $t+2<k$. 
Note that $U_1,\dots,U_{k-2}$ have duplication number $k-\ell-2$, so the expectation of the contribution of a set $U_i$ for $i\in[k-2]$ to the duplication number is less than 1. Now we choose $k-t-2$ sets $U_i$ from $\{U_1,\dots,U_{k-2}\}$, say $U_{t+1},U_{t+2}, \dots ,U_{k-2}$, with duplication number as large as possible. Then the duplication number of $U_{t+1},\dots,U_{k-2}$ is at most $t$ less than that of $U_1,\dots,U_{k-2}$. Recall the duplication number of $\mathscr{U}_{k-2}$ is at least $t+2$.
Therefore, the duplication number of $U_{t+1},U_{t+2}, \dots ,U_{t+k}$ is at least
 $(k-\ell-2)-t+(t+2)= k-\ell.$
Again a contradiction. The claimed assertion follows.

Assume that $|U_{k-1}\cap (\cup_{k-1\le j\neq i\le t+k}U_j)|\le 1$ without loss of generality. 
Since any $t+1$ elements of $\mathscr{U}$ have duplication number at most $b$,
we obtain $\sum_{i=k-1}^{t+k}| U_{i}|-|\cup_{k-1\le i\le t+k}U_i|\le b+1$.
\hfill
$\square$

\medskip
From the above proposition, we obtain that 
\begin{align*}
\sum_{i=1}^{t+k}| U_{i}|-|\cup_{1\le i\le t+k}U_i|
&=\left(\sum_{i=1}^{k-2}| U_{i}|-|\cup_{1\le i\le k-2}U_i|\right)+\left(\sum_{i=k-1}^{t+k}| U_{i}|-|\cup_{k-1\le i\le t+k}U_i|\right)
\\&\leq (k-\ell-2)+(b+1)
\\&= k-\ell-1+b.
\end{align*} 
So the claim holds for this case.

\medskip
(ii) Any $k-2$ elements of $\mathscr{U}$ have duplication number at most $k-\ell-3$. Similarly, there are two subcases. 

We repeat the above procedure at most $k-\ell-1$ times.
At each of these steps, we can inductively deduce that $\sum_{i=1}^{t+k}| U_{i}|-|\cup_{1\le i\le t+k}U_i|\le k-\ell-1+b.$ Indeed,
suppose at the $d$th step, any $k-d$  elements of $\mathscr{U}$ have duplication number at most $k-\ell-d-1$. We first assume that there exist $k-d-1$ sets $U_i$, say $U_1,\dots,U_{k-d-1}$, have duplication number $k-\ell-d-1$. Then by a similar argument as in Proposition \ref{indt}, we can  obtain by induction that
\begin{align*}
\sum_{i=k-d}^{t+k}| U_{i}|-\left|\bigcup_{k-d\le i\le t+k}U_i\right|
\le\left(\sum_{i=k-d+1}^{t+k}| U_{i}|-\left|\bigcup_{k-d+1\le i\le t+k}U_i\right|\right)+1
 \leq (b+d-1)+1=b+d
\end{align*} 
  It follows that
$
\sum_{i=1}^{t+k}| U_{i}|-|\cup_{1\le i\le t+k}U_i|
\leq (k-\ell-d-1)+(b+d)
= k-\ell-1+b.
$
 At the final step, i.e., the $(k-\ell-1)$th step, we may assume that any $\ell+1$ sets $U_i$ have duplication number equals to $0$. This implies that $\sum_{i=1}^{t+k}| U_{i}|-|\cup_{1\le i\le t+k}U_i|=0$. So the claim is proved.\hfill$\Box$
 
 \medskip
In the following, we aim to compute the maximum value of the  duplication number of $\mathscr{U}$. From the proof of the above claim, we may assume that there are $k-1$ elements of $\mathscr{U}$, say $U_{1},\dots,U_{k-1}$, whose  duplication number is
 $k-\ell-1$, and $\sum_{i=1}^{t+k}| U_{i}|-|V|$ achieves the maximum $k-\ell-1+b$. We will estimate the duplication number $b$ of the remaining $t+1$ sets. Let $\mathscr{A}=\{U_{1},\dots,U_{k-1}\}$ and $\mathscr{B}=\{U_{k},\dots,U_{t+k}\}$, and let $A=\cup_{i=1}^{k-1} U_{i}$ and $B=V\setminus A=\cup_{i=k}^{t+k} U_{i}$. 
 

Recall $U_i=V_i\cup W_i$ and note that only vertices from $W_i$ may contribute the duplication number. Note also that the subgraph of $R$ induced by $U_i\cup U_j$ for some distinct $i$ and $j$ would be connected if $W_i\cap W_j\neq\emptyset$.
Since the duplication number of any $k$ elements of $\mathscr{U}$ is at most $k-\ell-1=(k-1)-\ell$, we obtain that the subgraph induced by $A$ in $R$ has at least $\ell$ connected components.

Let $A_1$ be a connected component of $R[A]$. It consists of some sets from $\mathscr{A}$, and we denote this collection of sets by $\mathscr{A}_1$. Without loss of generality, we assume that $\mathscr{A}_1$ achieves the minimum duplication rate among these collections corresponding to the connected  components of $R[A]$. 
\begin{proposition}\label{dr-A1}
The duplication rate of the sets $U_i$ corresponding to any connected  component in $B$ should be less than or equal to that of $\mathscr{A}_1$.     
\end{proposition}
{\bf Proof.}
Suppose to the contrary that there is a connected  component $B_1$ in $B$ that contradicts the assumption. 
Let $\mathscr{B}_1$ be the collection of sets $U_i$ from $B_1$.
Suppose that $\mathscr{A}_1$ consists of $h_1$ sets with duplication rate $\rho_1$, and $\mathscr{B}_1$ consists of $h_2$ sets with duplication rate $\rho_2$. Then $\rho_2>\rho_1$. We claim that the duplication number of $\mathscr{A}_1$ is $h_1-1$. From the fact that $\rho_1<1$ as the duplication rate of  $\mathscr{A}$ is $\frac{k-\ell-1}{k-1}<1$, we obtain that the duplication number of $\mathscr{A}_1$ is at most $h_1-1$. On the other hand, since the subgraph induced by $A_1$ is connected, the duplication number of $\mathscr{A}_1$ is at least $h_1-1$. Thus the duplication number of $\mathscr{A}_1$ equals $h_1-1$.

\ignore{\begin{figure}[t]
\centerline{
\includegraphics[scale=0.67]{h1,2.pdf}}

\vspace{-1cm}
\caption{The case $h_1\ge h_2$}
\label{fg8}
\end{figure}
}
If $h_1\ge h_2$, then $\mathscr{B}_1$ has duplication number at least $h_2$ as $\rho_2>\rho_1=\frac{h_1-1}{h_1}$.
Thus we can take $h_1-(h_2-1)$ sets $U_i$ from $\mathscr{A}_1$ with maximum duplication rate, which have duplication number at least $h_1-h_2$ since the subgraph induced by $A_1$ is connected (That is, the duplication number is reduced by a maximum of $h_2-1$ compared to $\mathscr{A}_1$). Recall $\mathscr{A}$ has duplication number $k-\ell-1$. Thus, these $h_1-(h_2-1)$ sets $U_i$ from $\mathscr{A}_1$, together with the sets $U_i$ from $\mathscr{A}\setminus\mathscr{A}_1$ and the sets $U_i$ from  $\mathscr{B}_1$, yield $(k-1)-(h_2-1)+h_2=k$ elements of $\mathscr{U}$, whose duplication number is at least 
$(k-\ell-1)-(h_2-1)+h_2= k-\ell,$  a contradiction. 

If $h_1< h_2$, then we may take $h_1+1$ sets $U_i$ from $\mathscr{B}_1$ with maximum duplication rate and so they have duplication number at least $h_1$ as the subgraph induced by $B_1$ is connected. These sets, together with all sets $U_i$ from $\mathscr{A}\setminus\mathscr{A}_1$, yield $k$ elements of $\mathscr{U}$ with the duplication number at least $(k-\ell-1)-(h_1-1)+h_1=k-\ell$,  again a contradiction. \hfill$\Box$

\medskip 

To get the maximum value of $\sum_{i=1}^{t+k}| U_{i}|-|\cup_{1\le i\le t+k}U_{i}|$, it is necessary to achieve the maximum number of duplication number in $\mathscr{B}$, i.e., the duplication rate of the corresponding sets $U_i$ of each connected component in $B$ should reach the maximum. To achieve this, from Proposition \ref{dr-A1}, $A$ should have the minimum number of connected components and the duplication rates of the corresponding sets $U_i$ of all connected  components in $A$ are as close to the average value as possible. Since $\mathscr{A}$ has duplication number $k-\ell-1$, and $k-1=\mu \ell+\alpha$ where $0\le \alpha<\ell$ from the assumption, we may assume that $A$ contains exactly $\ell$ connected components, among which there are $\alpha$ components each contains $\mu+1$ sets $U_i$ for $i\in[t+k]$ with duplication number $\mu$ and $\ell-\alpha$ components each contains $\mu$ sets $U_i$ for $i\in[t+k]$ with duplication number $\mu-1$. In the above construction, each of the connected components in $B$ may reach the maximum duplication rate $(\mu-1)/\mu$ from Proposition \ref{dr-A1} by noting $(\mu-1)/\mu<\mu/(\mu+1)$. Thus the $t+k$ sets $U_i$ in $\mathscr{U}$ reaches the maximum duplication number, and so we obtain that
$$\sum_{i=1}^{t+k}| U_{i}|-|\cup_{1\le i\le t+k}U_{i}|\leq\alpha\mu+\tfrac{\mu-1}{\mu}(t+k-\alpha),$$ 
which implies that $|V|\geq \sum_{i=1}^{t+k}| U_i|-\tfrac{\mu-1}{\mu}(t+k-\alpha)-\alpha\mu$ as claimed. \hfill $\square$

\section{Proof of Theorem \ref{main4}}\label{ub}
 
 The lower bounds of Theorem \ref{main4} have been established in Section \ref{lb}, so we only need to prove the upper bounds accordingly.
In the following, for each case, we consider a red-blue edge coloring of \(K_N\) for some suitable $N$ defined on the vertex set \(V\), and we consistently use $R$ and $B$ to denote the graphs induced by all red edges and blue edges respectively.
Let
\[
p=\left\lfloor (n-1)/t \right\rfloor, \;\;\text{and so} \;\;n-1-t< tp\le n-1.
\]
For any $t\ge 2$ and $k\ge3$, let $s\ge k$, and let $K$ and $m$ be defined as in Section \ref{prop-ex}.

\subsection{The upper bound for Theorem \ref{main4}(i)}

From Theorem \ref{hllnp}, we know Theorem \ref{main4}(i) holds for $k=2$, we now assume that $k\ge3$ and Theorem \ref{main4}(i) holds for all smaller $k$. We will show that Theorem \ref{main4}(i) also holds for $k$.

\medskip\noindent
{\bf Case 1.} \ $g_{k-1}=(t+k-2)(m-1)+1.$

\medskip
For this case, $n+(k-1)p+k-1\leq(t+k-2)(m-1)+1$. 
Note that $p+1<m-1$, so we obtain $n+kp+k<(t+k-1)(m-1)+1.$ Let \(N=(t+k-1)(m-1)+1 \). We will show $g_k\le N$. Consider a red-blue edge coloring of \(K_N\) on vertex set \(V\). On the contrary, suppose that the red graph $R$ contains no $C_m$ and the blue graph $B$ contains no $B_n^{(k)}$.

If there exists a vertex \(v\in V\) with \(d_B(v)\geq(t+k-2)(m-1)+1=g_{k-1}\), then by the inductive hypothesis we can find either a red \(C_m\) or a blue \(B_n^{(k-1)}\) in the blue neighborhood of $v$. In the latter case, such a \(B_n^{(k-1)}\) together with $v$ gives a blue \(B_n^{(k)}\).

So we may assume $\delta(R)\geq m - 1$. If \(R\) is 2-connected, then since $\delta(R)\geq m - 1$, $N>2m-2+1$, by Lemma \ref{VZ}, $ec(R)\geq 2m-2$. Thus $R$ contains $C_m$ by Lemma \ref{GHS}. 

Now we assume that $R$ has a cut-vertex. Note that
$
    \delta(R)\geq m - 1\geq\frac{(t+k-1)(m-1)+1}{t+k}+(t+k),  
$
provided $m\ge(t+k)^2$.
By Lemma \ref{Luzack} there is an $s \leq t+k-1$ and a set $W\subset V$ with $|W| \leq s-1$ such that $R-W$ is a union of vertex-disjoint 2-connected subgraphs $V_1,\ldots,V_s$. Let $U_i$ be a subset such that $U_i=V_i\cup W_i$, where $W_i=\{w\in W : d_R(w,V_i)\ge 2s\}$. Since $\delta(R)\ge m-1\ge m/2+K/2+s$, by Lemma \ref{basic}(i)-(ii) we know that $V=\cup_{i=1}^s U_i$ and $|U_i| \leq m-1$. However, then 
$
|V|=|\cup_{i=1}^s U_i|\leq(t+k-1)(m-1)<N,
$
which leads to a contradiction.

\vspace{0.2cm}
\noindent
{\bf Case 2.} \ $g_{k-1}=n+(k-1)p+(k-1).$

\vspace{0.2cm}

We first show that $g_k\le n+kp+k.$
Let $N=n+kp+k$ and consider a red-blue edge coloring of \(K_N\) on the vertex set \(V\). If there exists a vertex \(v\in V\) with \(d_B(v)\geq n+(k-1)p+(k-1)=g_{k-1}\), then we are done similarly as above.
So we assume $\delta(R)\geq p+1$.
If $t=2$, by the induction hypothesis, $n+(k-1)p+(k-1)\ge k(m-1)+1$, implying $n\ge \tfrac{2k}{k+1}(m-1)-\tfrac{k-3}{k+1} \ge  3(m-1)/2$ as $k\ge 3$. Thus $2\delta(R)\ge 2p+2\ge n \ge 3(m-1)/2.$ If $t\ge 3$, then as $n-1\ge (t-1)(m-1)$, we know $2\delta(R)\ge 2p+2>2(n-1)/t \ge 2(t-1)(m-1)/t \ge 4(m-1)/3.$ Therefore we conclude that $2\delta(R)\ge 4(m-1)/3$.

If $R$ is 2-connected, then by a similar argument as above, $R$ contains $C_m$. 
Thus we assume that $R$ has a cut-vertex.
Note that
$
    \delta(R)\geq p+1\geq\frac{N}{t+k+1}+(t+k+1),
$
provided $n\ge2t(t+k+1)^2$. By Lemma \ref{Luzack}, there is an $s \leq (t+k)$ and a set $W\subset V$ with $|W| \leq s-1$ such that $R-W$ is a union of vertex-disjoint 2-connected subgraphs $V_1,\ldots,V_s$.
For each $1 \leq i \leq s$, we define $U_i$ as before. Since $\delta(R)\ge p+1\ge 2(m-1)/3>m/2+K/2+s$, by Lemma \ref{basic}(i)-(ii) we know that $V=\cup_{i=1}^s U_i$ and $|U_i| \leq m-1$.

If $s\leq t+k-1$, then
$|V|=|\cup_{i=1}^s U_i|\leq(k+t-1)(m-1)<(k+t-1)(m-1)+1\le N,
$
a contradiction.
If $s=t+k$, as $\delta(R)\geq p+1$, then $|U_i|\geq p+2$. {By Lemma \ref{basic}(iii) and $\delta(R)\ge p+1\ge 2(m-1)/3>m/2+K/2+s$,} we obtain
\begin{align*}
|V|\geq\sum_{i=1}^s| U_i|-(s-1)
         \geq(t+k)(p+2)-(t+k-1)\ge n+kp+k+1>N,
\end{align*}
where the third inequality holds since $t(p+1)\ge n$, a contradiction.

Now we assume $(t+k-1)(m-1)+1\geq n+kp+k$ and we will show $g_k\le(t+k-1)(m-1)+1.$
Let $N=(t+k-1)(m-1)+1$, and consider an arbitrary red-blue edge coloring of $K_N$ on the vertex set $V$. 
We rephrase $N$ as 
$
    N=n+kp+k+a-1,
$
where $a\geq 1$. If \(d_B(v)\geq n+(k-1
) p+(k-1)=g_{k-1}\) for some $v\in V$, then we are done similarly as above. So we assume $\delta(R)\geq p+a$. 

If $R$ is 2-connected, then by a similar argument as above we can obtain a red $C_m$. 
So we may assume that $R$ contains a cut-vertex. 
Note that 
$
 \delta(R)\geq p+a\geq\frac{N}{t+k+1}+(t+k+1).   
$
By Lemma \ref{Luzack}, there is an $s \leq (t+k)$ and a set $W$ with $|W| \leq s-1$ such that $R-W$ is a union of vertex-disjoint 2-connected subgraphs $V_1,\ldots,V_s$.
For each $1 \leq i \leq s$, we define $U_i$ as before. As $\delta(R)\ge p+a\ge m/2+K/2+s$, by Lemma \ref{basic}(i)-(ii) we know that $V=\cup_{i=1}^s U_i$ and $|U_i| \leq m-1$. 
If $s\leq t+k-1$, then $|V|=|\cup_{i=1}^s U_i|\leq(k+t-1)(m-1)<N,$
a contradiction. 
If $s=t+k$, then $|U_i|\geq p+a+1$ since $\delta(R)\geq p+a$. Thus, by Lemma \ref{basic}(iii),
\begin{align*}
|V|\geq\sum_{i=1}^s| U_i|-(s-1)
         &\geq(t+k)(p+a+1)-(t+k-1)\\
             &\ge  n+kp+ak+(a-1)t+1\\
       &> N,
\end{align*}
where the third inequality holds since $t(p+1)\ge n$, 
a contradiction.

\subsection{The upper bound for Theorem \ref{main4}(ii)}
The proof proceeds by induction on $k\ge2.$ The base case follows directly from Theorem \ref{hllnp}. For the inductive step, assume $k\ge3$ and that the assertion holds for all smaller values of $k$. We will show that the assertion also holds for $k$.

Suppose $p = m-\sigma$, where $ \lceil \frac{k}{2}\rceil+1\le \sigma\le k$. 
Note that \begin{align*}
    n+(k-1)p=(t+k-1)p+n-tp\ge(t+k-1)(m-\sigma)
    >(t+k-2)(m-1)+1,
\end{align*}
provided $m\ge k(t+k)$.
If $p=m-k=m-(k-1)-1$, then it follows by Theorem \ref{main4}(i) that $g_{k-1}=n+(k-1)p+k-1$.
If $p = m-\sigma$, where $\lceil \frac{k}{2}\rceil+1\le \sigma\le k-1 $, by induction, we assume that
$g_{k-1}=n+(k-1)p+\sigma-1$.
We will show  $g_{k}\le n+kp+\sigma-1$ for $p = m-\sigma$, where $ \lceil \frac{k}{2}\rceil+1\le \sigma\le k$.

Let $N= n+kp+\sigma-1$, and consider a  red-blue edge coloring of \(K_N\) on the vertex set \(V\). If there is a vertex $v$ with \(d_B(v)\geq n+(k-1)p+\sigma-1=g_{k-1}\), then we are done similarly as above.
So we may assume $\delta(R)\geq p$. If \(R\) is 2-connected, by a similar argument, $R$ contains $C_m$.

Thus we assume that $R$ contains a cut-vertex.
Note that $\delta(R)\geq p\geq\frac{N}{t+k+1}+(t+k+1)$. By Lemma \ref{Luzack}, there is an $s \leq (t+k)$ and a set $W$ with $|W| \leq s-1$ such that $R-W$ is a union of vertex-disjoint 2-connected subgraphs $V_1, \ldots, V_s$.
For each $1 \leq i \leq s$, we define $U_i$ as before. Since $\delta(R)\ge m-k\ge m/2+K/2+s$, by Lemma \ref{basic}(i)-(ii) we know that $V=\cup_{i=1}^s U_i$ and  $|U_i| \leq m-1$. As $\delta(R)\ge p$, we have that $$p+1\le |U_i|\le p+\sigma-1.$$
If $s\leq t+k-1$, then
$
|V|\leq\sum_{i=1}^s| U_i|
         \le (t+k-1)(p+\sigma-1)\le n+kp-p+(t+k-1)(\sigma-1)<N,
$
where the last inequality holds provided $n\ge tk(t+k)$, which leads to a contradiction.
Thus we may assume $s=t+k$. 

\begin{claim}\label{min-s-v}
The size of $V$ achieves the minimum when $|U_i|=p+1$ for all $1\le i\le t+k$.
\end{claim}
{\bf Proof.}
 For $|U_i|\ge p+2$, and for $|U_j|=p+y$, $1\le y\le \sigma-1$, $1\le j\not=i\le t+k$, set $|U_i\cap U_j|=x$. If $x\le y-1$, then we can choose $U_j'\subseteq U_j$ of $p+y-x$ vertices that do not intersect with $U_i$. And if $x\ge y$, we can choose $U_j'\subseteq U_j$ of $p+1$ vertices that intersects $U_i$ with $x+1-y\ge1$ vertices. 
  For each of these cases, if we substitute $U_j$ by $U_j'$, it will not change the size of $V$.
  
Let $\lambda$ be the number of $U_j'$'s that intersects $U_i$ with at least one vertex. If $\lambda\ge k-1$, then there are $U_{j_1},\dots,U_{j_{k-1}}$ corresponding to these $U_j'$'s such that 
\begin{align*}
 |U_i\cup U_{j_1} \cup \cdots \cup U_{j_{k-1}}| \le (p+\sigma-1)+ (k-1)p=kp+\sigma-1,
\end{align*}
implying that $|V\setminus(U_i\cup U_{j_1} \cup \cdots \cup U_{j_{k-1}})|\ge n$, which contradicts Lemma \ref{c5}. 

So we assume $\lambda< k-1$.
Let $\tau=|\cup_{j\neq i} (U_j'\cap U_i)|.$
If $\tau\le|U_i|-p-1$, then we can choose $U_i'\subseteq U_i$ with $|U_i|-\tau\ge p+1$ vertices which is disjoint with all other $U_j'$'s. Thus the size of $V$ achieves the minimum when $|U_i'|=p+1$.
If $\tau>|U_i|-p-1$, then we can choose $U_i'\subseteq U_i$ of $p+1$ vertices that intersects $\tau+(p+1)-|U_i|$ vertices in total with all other $U_j'$'s. In particular, these $\tau+(p+1)-|U_i|$ vertices are contained in $V\setminus U_i'$. So if we replace $U_i$ by $U_i'$, it will not effect the size of $V$ for this case.  
Therefore, we conclude that the size of $V$ achieves the minimum when $|U_i|=p+1$ for all $1\le i\le t+k$. \hfill$\Box$

\medskip

By Lemma \ref{c5}, the size of any $k$ sets $U_i$ for $1\le i\le t+k$ is at least $N-(n-1)=kp+\sigma$.
Now we choose $k$ sets $U_i$ for $1\le i\le t+k$ that have maximal duplication number, which is at most $k-\sigma\le \lfloor\tfrac{k}{2}\rfloor-1$ because $ \lceil \frac{k}{2}\rceil+1\le \sigma\le k $. Among these $k$ sets $U_i$, there are at least two disjoint sets of them that are disjoint from all the others. Moreover, for all other $U_j$'s, they are surely disjoint from each other and also disjoint from these $k$ sets $U_i$. Otherwise, we can substitute two disjoint sets from these $k$ sets $U_i$'s by two intersecting sets out of them, yielding a new cluster of $k$ sets with larger duplication number, a contradiction.
Therefore, 
 \begin{align*}
|V|\geq(t+k)(p+1)-(k-\sigma)>n-1+kp+\sigma=N,
\end{align*}
which leads to a contradiction.

\subsection{The upper bound for Theorem \ref{main4}(iii)}

Suppose $p = m-\sigma$, where $ 2\le \sigma\le \lceil \frac{k}{2}\rceil$, and $g_{k-1}=n+(k-1)p+\ell$ for some $\sigma-1 \le \ell \le \left\lceil \tfrac{k-1}{2} \right\rceil$, and $k-1 = \mu_\ell \ell + \alpha_\ell$ where $0\le\alpha_\ell <\ell$, we will show that
\[
g_{k}\le \begin{cases}
n+kp+\ell+1 &  \textrm{if $r_{k}\le r$,}\\
n+kp+\ell &  \textrm{if  $r_{k}>r$.}
\end{cases}
\]
where $r_{k}=tp+t+k-\ell-n$, and $r=\tfrac{\mu_\ell-1}{\mu_\ell} (t+k-\alpha_\ell)+\alpha_\ell \mu_\ell$.

\vspace{0.2cm}
\noindent
{\bf Case 1.} \ $r_{k}\le r$.

\vspace{0.2cm}
Let $N=n+kp+\ell+1$, and consider a red-blue edge coloring of \(K_N\) on \(V\). If there is a vertex $v$ with \(d_B(v)\geq n+(k-1)p+\ell=g_{k-1}\), we are done similarly.
So we assume  $\delta(R)\geq p+1$. 
If \(R\) is 2-connected, then by a similar argument as above, $R$ contains $C_m$ by Lemma \ref{GHS}. 

Now we suppose that $R$ contains a cut-vertex. Note that $\delta(R)\geq p+1\geq\frac{N}{t+k+1}+(t+k+1)$. By Lemma \ref{Luzack}, there is an $s \leq (t+k)$ and a set $W$ with $|W| \leq s-1$ such that $R-W$ is a union of vertex-disjoint 2-connected subgraphs $V_1, \ldots, V_s$.
For each $1 \leq i \leq s$, we define $U_i$ as before. Since $\delta(R)\ge m-k\ge m/2+K/2+s$, by Lemma \ref{basic}(i)-(ii) we know that $V=\cup_{i=1}^s U_i$ and $|U_i| \leq m-1$. As $\delta(R)\ge p+1$, then $p+2\le|U_i|\le p+\sigma-1$. 

If $s\le t+k-1$, then by a similar calculation as in Theorem \ref{main4}(ii) we obtain that
\begin{align*}
|V|&\leq\sum_{i=1}^s| U_i|\le (t+k-1)(p+\sigma-1)<N,
\end{align*}
provided $n\ge tk(t+k)$,  a contradiction.

If $s= t+k$, then by Lemma \ref{basic}(iii) with $\delta(R)\ge m-k\ge m/2+K/2+s$, we obtain that
\begin{align*}
|V|\geq\sum_{i=1}^s| U_i|-(s-1)
         \geq(t+k)(p+2)-(t+k-1)\ge n+kp+k+1>N,
\end{align*}
where the third inequality holds since $t(p+1)\ge n$, a contradiction.

\vspace{0.2cm}
\noindent
{\bf Case 2.} \ $r_{k}>r$.

\vspace{0.2cm}
Let $N=n+kp+\ell$, and consider a  red-blue edge coloring of \(K_N\) on \(V\). If there is a vertex $v$ that \(d_B(v)\geq n+(k-1)p+\ell= g_{k-1}\), we are done as before.
So we assume $\delta(R)\geq p$. 

If \(R\) is 2-connected, then by a similar argument as above, we can obtain a $C_m$ in $R$. 
Thus we assume that $R$ contains a cut-vertex.
Note that $\delta(R)\geq p\geq\frac{N}{t+k+1}+(t+k+1)$. It follows by Lemma \ref{Luzack} that there is an $s \leq (t+k)$ and a set $W$ with $|W| \leq s-1$ such that $R-W$ is a union of vertex-disjoint 2-connected subgraphs $V_1, \ldots, V_s$.
For each $1 \leq i \leq s$, we define $U_i$ as before. Since $\delta(R)\ge m-k\ge m/2+K/2+s$, by Lemma \ref{basic}(i)-(ii) we know that $V=\cup_{i=1}^s U_i$ and $|U_i| \leq m-1$. As $\delta(R)\ge p$, we have $p+1\le|U_i|\le p+\sigma-1$. 

If $s\leq t+k-1$, then
$
|V|\leq\sum_{i=1}^s| U_i|\le (t+k-1)(p+\sigma-1)<N,
$
a contradiction.
So we assume $s=t+k$. Since $\delta(R)\ge m-k\geq 2s^2+s$ and $B$ contains no $B_n^{(k)}$. By Lemma \ref{c5}, for any distinct $i_1,i_2,\dots,i_{k}\in [s]$, we have $|\cup_{1\le j\le k} U_{i_j}|\geq |V|-(n-1)= kp+\ell+1$. By a similar argument as Claim \ref{min-s-v}, we obtain that the size of $V$ achieves the minimum when $|U_i|=p+1$ for all $1\le i\le t+k$.
Therefore, by Lemma \ref{c6},
\begin{align*}
|V|&\geq(t+k)(p+1)-(\tfrac{\mu_\ell-1}{\mu_\ell} (t+k-\alpha_\ell)+\alpha_\ell \mu_\ell)\\
         &>(t+k)(p+1)-r_k\\
       &=n+kp+\ell\\
    &=N,
\end{align*}
where the first equality holds since $r_k=tp+t+k-\ell-n$.
The final contradiction completes the proof. \hfill$\Box$

\ignore{\medskip

{\em Remark.} In the proof of Case 1 for the lower bound of Theorem \ref{main4}(iii), we know that $r_k=tp+t+k-\ell-n$ represents the duplication number of the vertices in all $K_{p+1}$'s in the construction of $\Gamma_k$. Meanwhile, Lemma \ref{c6} states that $r=\tfrac{\mu_\ell-1}{\mu_\ell} (t+k-\alpha_\ell)+\alpha_\ell \mu_\ell$ is the maximum duplication number of the vertices in all  $K_{p+1}$'s in a $(C_m,B_{n}^{(k)})$-Ramsey graph. 
}

\section{Concluding remarks}

\noindent

{$\bullet$} In the proof of Case 1 for the lower bound of Theorem \ref{main4}(iii), $r_k=tp+t+k-\ell-n$ represents the duplication number of the vertices in all $K_{p+1}$'s in the construction of $\Gamma_k$. Meanwhile, Lemma \ref{c6} states that $r=\tfrac{\mu_\ell-1}{\mu_\ell} (t+k-\alpha_\ell)+\alpha_\ell \mu_\ell$ is the maximum duplication number of the vertices in all  $K_{p+1}$'s in a $(C_m,B_{n}^{(k)})$-Ramsey graph. 

\medskip
$\bullet$ Let us point out that
Lemma \ref{GHS} by Gould, Haxell, and Scott \cite{GHS02} can be slightly improved as follows \cite{scott}:

\medskip
{\em Let $c> 0$, $K=75\cdot 10^4/c^5$, and let $G$ be a graph with $n\geq {18K}/{c^2}$
vertices and minimum degree at least $c n$. Then $G$ contains a cycle of length $t$ for
every even integer $t\in [4,ec(G)-K]$.}

\medskip
Indeed, in the proof of \cite[Theorem 1]{GHS02}: ``$\cdots$Otherwise at least $2c^2n/9 - 2K_0$ common neighbours of $z$ and $z_0$ fall
onto $P_i$. Since $P_i$ has length less than $n$, and $n\ge90K_0c^{-4}$, there is an interval $I$
in $P_i$ of length at most $15c^{-2}$ that is disjoint from $S$ and contains three common
neighbours and therefore two that are an even distance apart in $I$, say $y_1$ and $y_2$$\cdots$'', here the condition ``$n\ge90K_0c^{-4}$'' is not required and can be relaxed to ``$n\ge90K_0c^{-2}$'', say. Moreover, this condition is not needed anywhere else. 

\medskip
From the above observation, we can thus slightly improve Theorem \ref{Thm:alpz}, Theorem \ref{hllnp}, and Theorem \ref{main4} to obtain the exact value of $R(C_m,B_{n}^{(k)})$ for each fixed $k\ge1$ and large even $m\ge \Omega(n^{7/8})$. However,   these three theorems  do not hold if $m\le O(\log n/\log\log n)$ from (\ref{bou-g}) as was observed in \cite{ALPZ-Arxiv}.
Finding the minimal size of $m$ for which this holds remains an interesting problem. 

\end{spacing}

\end{document}